\documentclass[12pt,a4paper]{amsart}
\usepackage{amsmath}
\usepackage{pstricks,pst-node}

\newtheorem{defn}{Definition}[section]

\newtheorem{example}{Example}[section]

\newtheorem{prop}{Proposition}[section]
\newtheorem{lemma}{Lemma}[section]

\newtheorem{theorem}{Theorem}[section]

\begin{document}
\title{Enumeration of non-positive planar trivalent graphs}
\author{Bruce W. Westbury}
\address{Mathematics Institute\\
University of Warwick\\
Coventry CV4 7AL}
\email{bww@maths.warwick.ac.uk}
\date{\today}
\begin{abstract} In this paper we construct inverse bijections
between two sequences of finite sets. One sequence is defined
by planar diagrams and the other by lattice walks. In
\cite{MR1403861} it is shown that the number of elements in
these two sets are equal. This problem and the methods we use
are motivated by the representation theory of the exceptional
simple Lie algebra $G_2$. However in this account we have
emphasised the combinatorics.
\end{abstract}
\maketitle

\section{Introduction}
The aim of this paper is to give an enumeration of non-positive
planar trivalent graphs. Here the graph is embedded in a disk
and the number of boundary points is specified. This differs from
other enumerations of planar trivalent graphs such as
\cite{MR0130841} or \cite{MR1980342} in that neither the number of
vertices not the number of edges is specified. So if there were no
further conditions then the number of graphs would be infinite.
The extra condition that is imposed is that the graph is non-positive.
This means that there is no internal face with less than six edges.

This problem arose in the work of Greg Kuperberg on the representation
theory of the exceptional simple Lie algebra $G_2$ in 
\cite{MR1403861}. In particular, one of the results of this paper
is that the two sets we consider have the same numbers of elements.
This is proved by showing that both sets are a basis of the same
vector space.

Here, we give a bijective proof of this result. A bijective
proof of the analogous result for $A_2$ is given in 
\cite{MR1680395}. Moreover, the map from diagrams to words given below
is constructed by the same method as the analogous map in this reference.
However the construction of the inverse map that we give here is new.
This construction is based on a diagram model for the crystal graph.

The main result of this paper is the construction of inverse bijections
between two sequences of finite sets. The two sets and the bijections
are constructed by combinatorial methods and in writing this paper we
have emphasised the combinatorial aspects. However both the original
problem and the construction of the bijections are motivated by 
two combinatorial methods in representation theory, namely Littelman paths
and crystal graphs. So here we put the constructions we have used into
context.

In this paper we are studying the invariant theory of the seven
dimensional fundamental representation of the exceptional simple
Lie algebra of type $G_2$. The Lie algebra can be constructed as the
derivation algebra of the octonions and the representation can then
be taken to be the imaginary octonions. The weights and the dominant
weights below are the usual weights and dominant weights of this Lie
algebra. This representation has the (rare) property that all weight
spaces are one dimensional. The set $S$ in \eqref{steps} consists of the
seven weights
of this representation; moreover Figure \ref{dst} (with the centre point
of weight $(0,0)$) is the weight diagram of this representation.

The two finite sets we consider here each describe a basis for
the vector space of invariants in the $n$-th tensor power of this
representation. One basis is given by dominant Littelman paths.
The Littelman paths in this example consist of the six straight line
paths which are the edges emanating from the origin in Figure \ref{dst}
together with one path of weight zero. Then the lattice walks are
exactly the dominant paths, that is paths which are obtained by
concatenations of these seven paths which are dominant at all times.
The second method for constructing a basis of invariant tensors
is to use crystal graphs. The crystal graph in this example is
a directed graph with edges labelled by the simple roots
where the vertices are the set $S$.
The crystal graph is obtained from the crystal base.
The crystal base and the crystal graph for this representation are
given in \cite{MR1265857}. In this example the vertices are 
parametrised by the weights and two weights are connected by an
edge labelled by a simple root if the difference is the simple root.
In this paper we do not make any use
of the edges in the crystal graph and so they have been omitted.
However the vertices of a crystal graph are also labelled by two
dominant weights which we denote $H$ and $D$ where the weight is
$H-D$. There is a tensor
product construction for crystal graphs and this includes a rule
for the labels $H$ and $D$. The labelling in Figure \ref{dw} is
designed to reproduce this rule.

\section{Lattice walks}
In this section we will give the description of the lattice walks
that will be related to the diagrams. We will use the terminology
of weights which comes from representation theory.

A weight $\lambda$ is an ordered pair of integers. The set of
weights is an abelian group under addition and is partially ordered.
The partial order is given by saying $(a_1,b_1)\ge(a_2,b_2)$ if
$a_1\ge a_2$ and $b_1\ge b_2$. A weight $\lambda$ such that 
$\lambda\ge 0$ is called a dominant weight.

Define the set $S$ to be the following set of seven weights
\begin{equation}\label{steps}
S = \{ (-2,1),(0,-1),(-1,1),(0,0),(1,0),(-1,1),(2,-1) \}
\end{equation}
Then the main definition of this section is
\begin{defn} A lattice walk of length $n$ is a sequence of
$n+1$ dominant weights $(\lambda_0,\lambda_1,\ldots ,\lambda_n)$
such that $(\lambda_i-\lambda_{i-1})\in S$ for $1\le i\le n$.
The sequence is also required to satisfy the additional condition
that, for $1\le i\le n$, if $\lambda_i=(a,0)$ for some $a\ge 0$
then $\lambda_{i-1}\ne\lambda_i$.
\end{defn}
If $(\lambda_0,\lambda_1,\ldots ,\lambda_n)$ is a lattice walk
then we say that the walk starts at $\lambda_0$ and ends at
$\lambda_n$. Denote the set of lattice walks of length $n$
which start at $\lambda$ and end at $\mu$ by $W(\lambda,n,\mu)$.
Also let the number of elements of $W(\lambda,n,\mu)$ be denoted
by $w(\lambda,n,\mu)$.

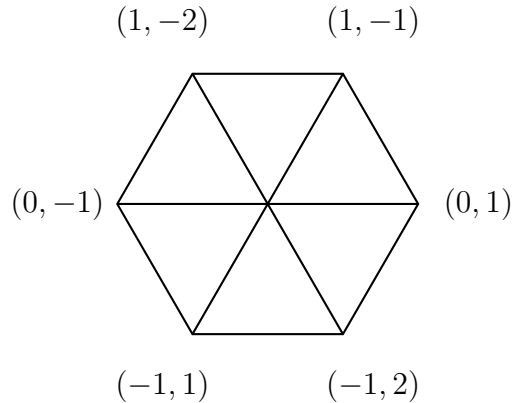
\begin{figure}
\begin{center}
\psset{yunit=0.866}
\begin{pspicture}(-3,-3)(3,3)
\pspolygon(2,0)(1,-2)(-1,-2)(-2,0)(-1,2)(1,2)
\psline(2,0)(-2,0)
\psline(1,2)(-1,-2)
\psline(1,-2)(-1,2)
\rput(2.8,0){$(0,1)$}
\rput(1.4,2.8){$(1,-1)$}
\rput(-1.4,2.8){$(1,-2)$}
\rput(1.4,-2.8){$(-1,2)$}
\rput(-1.4,-2.8){$(-1,1)$}
\rput(-2.8,0){$(0,-1)$}
\end{pspicture}
\psset{xunit=1cm,yunit=1cm}
\end{center}
\caption{Seven steps}\label{dst}
\end{figure}

Next we interpret a lattice walk as a walk in the triangular lattice.
First we identify the set of weights with the vertices of the
triangular lattice. This is an identification of abelian groups.
The identification we use identifies the non-zero weights in $S$
with the regular hexagon of nearest neighbours of the origin in the
triangular lattice. This is shown in Figure \ref{dst}.

\begin{example} The following diagram has $w(0,6,\mu)$ written at
position $\mu$.
\begin{equation*}
\begin{array}{lllllllllllll}
&  &  &  &  &  &  &  &  & 5 &  &  &  \\ 
&  &  &  &  &  & 65 &  & 40 &  & 9 &  &  \\ 
&  &  & 120 &  & 176 &  & 120 &  & 40 &  & 5 &  \\ 
35 &  & 120 &  & 180 &  & 145 &  & 65 &  & 15 &  & 1
\end{array}
\end{equation*}
\end{example}

It is clear that the sets $W(\lambda,n,\mu)$ can be enumerated.
Here we give numbers of the form $w(0,n,\mu)$.
Let $X$ and $Y$ be the following Laurent polynomials in the
indeterminates $x$ and $y$:
\begin{eqnarray}\label{lp}
X &=& 1+x+x^{-1}+x^{-1}y+xy^{-1}+x^{-2}y+x^2y^{-1} \\
Y &=& x - x^{-3}y^2 + x^{-6}y^3 - x^{-8}y^3 + x^{-8}y^2 \\
&& - x^{-7} + x^{-3}y^{-2} - xy^{-4} +x^4y^{-5}
-x^6y^{-5}+x^6y^{-4}-x^4y^{-2}
\end{eqnarray}

Then the numbers $w(0,n,\mu)$ can be computed from the observation
that if $\mu=(a,b)$ then $w(0,n,\mu)$ is the coefficient of
$x^{a}y^{b}$ in the Laurent polynomial $YX^n$.

The sequence $a(n)=w(0,n,0)$ which is the number of lattice walks
of length $n$ which start and end at the origin is of particular
interest. This sequence is given in \cite{oeis} as sequence
A059710. Equivalently $a(n)$ is the constant term in the
exapnsion of $YX^n$ where $X$ and $Y$ are the Laurent polynomials
in \eqref{lp}. It follows from this description that the sequence
$a(n)$ is P-recursive which means that it satisfies a finite
recurrence relation with polynomial coefficients. Alex Mihailovs
has proposed that this sequence satisfies the following recurrence
relation
\begin{multline*} (n+5)(n+6)a(n)= 2(n-1)(2n+5)a(n-1)
\\+(n-1)(19n+18)a(n-2)+14(n-1)(n-2)a(n-3)
\end{multline*}
together with the initial conditions $a(0)=1$, $a(1)=0$, $a(2)=1$.
The evidence for the proposal is that it gives the first thirty
terms of the sequence and gives the correct asymptotics.

\section{Diagrams}
The main definition of this section is the following:
\begin{defn}\label{md} A diagram with $n$
boundary points consists of a disc with $n$ marked points on the
boundary together with an embedded graph. The graph has $n$ vertices
of valence one which are identified with the marked boundary points
by the embedding and all other vertices of the graph have valency three.
\end{defn}
A region is a connected component of the complement of the image of
the graph in the disc.
\begin{defn} A diagram is non-positive if every 
region of the disc which is bounded by edges of the graph is
bounded by at least six edges of the graph.
\end{defn}

\begin{example} There are four non-positive diagrams with
four boundary points. These are
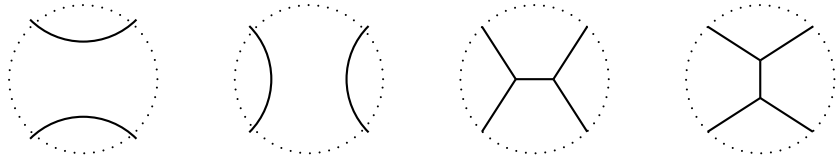
\begin{center}
\begin{pspicture}(-2,-1)(11,2)
\psset{linestyle=dotted}
\pscircle(0,0){1}
\pscircle(3,0){1}
\pscircle(6,0){1}
\pscircle(9,0){1}
\psset{linestyle=solid}
\psarc(0,1.5){1}{-135}{-45}
\psarc(0,-1.5){1}{45}{135}
\psarc(1.5,0){1}{-45}{45}
\psarc(4.5,0){1}{135}{225}
\psline(5.3,0.7)(5.75,0)
\psline(5.75,0)(6.25,0)
\psline(5.75,0)(5.3,-0.7)
\psline(6.7,0.7)(6.25,0)
\psline(6.25,0)(6.7,-0.7)
\psline(8.3,0.7)(9,0.25)
\psline(9,0.25)(9,-0.25)
\psline(9,-0.25)(8.3,-0.7)
\psline(9.7,0.7)(9,0.25)
\psline(9,-0.25)(9.7,-0.7)
\end{pspicture}
\end{center}
\end{example}

Although we will not make use of this result we first recall
the argument from \cite{MR1403861} which shows 
that if we specify the number of
boundary points then there are only finitely many non-positive
diagrams. This proof depends on the isoperimetric inequality
given in \cite{MR936419}.
\begin{prop}
Let $n\ge 0$. Then there are finitely
many non-positive diagrams with $n$ boundary points.
\end{prop}
\begin{proof} Consider a diagram whose graph is connected.
Then the dual graph gives a
triangulation of the disc. Take each triangle to be a Euclidean
equilateral triangle with edge length 1. Then this gives a
polyhedral metric on the disc. Now if the planar graph is non-positive
then this polyhedral metric has non-positive curvature. Hence
the isoperimetric inequality is satisfied. Each triangle has area
$\sqrt{3}/2$ and the length of the boundary is $n$. Hence the
isoperimetric inequality gives that there at most $n^2/(\pi\sqrt{3})$
triangles.
\end{proof}

Next we recall a crucial definition from \cite{MR1403861}.
\begin{defn}\label{cp} Assume we are given a diagram.
Let $A$ and $B$ be two boundary points which are not marked points.
Then a cut path from $A$ to $B$ is a path from $A$ to $B$ such that
each component of the intersection with the embedded graph is either
an isolated tranverse intersection point or else is an edge of the
graph.
\end{defn}
The diagrams for these two cases are:
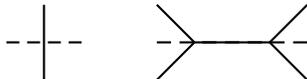
\begin{center}
\begin{pspicture}(0,-0.5)(4,0.5)
\psline[linestyle=dashed](0,0)(1,0)
\psline(0.5,0.5)(0.5,-0.5)
\psline[linestyle=dashed](2,0)(4,0)
\psline(2,0.5)(2.5,0)
\psline(2,-0.5)(2.5,0)
\psline(2.5,0)(3.5,0)
\psline(3.5,0)(4,0.5)
\psline(3.5,0)(4,-0.5)
\end{pspicture}
\end{center}

A cut path which crosses $a$ edges and contains $b$ edges is assigned
the weight $(a,b)$. The weight $(a,b)$ is called dominant if
$a\ge 0$ and $b\ge 0$. The weights are partially ordered by
$(a_1,b_1) < (a_2,b_2)$ if $a_1+2b_1 < a_2 + 2b_2$ or if
$a_1+2b_1 = a_2 + 2b_2$ and $b_1 < b_2$.
A cut path is minimal if there is no cut path with the same endpoints
and with lower weight.

%

\begin{defn} A triangular diagram is a non-positive diagram together with
three points $A$, $X$ and $Y$ which are in the boundary of the disc but
not marked points. Then we require
that the edges $AX$ and $AY$ are minimal cut paths.
\end{defn}
Our convention is that we draw a triangular diagram as a graph 
embedded in the triangle $AXY$ where $XY$ is a horizontal edge and
the vertex $A$ is below this edge.
\begin{defn}\label{rtd}
A triangular diagram is reducible if there is a point $B$ inside the
triangle (and not on the graph) such that there is a minimal cut path
from $A$ to $X$ which passes through $B$, a minimal cut path from
$A$ to $Y$ which also passes through $B$ and such  that the the triangular
diagram with vertices $B$, $X$ and $Y$ and edges given by these minimal
cut paths from $B$ to $X$ and $B$ to $Y$ is a proper subdiagram.
\end{defn}
A triangular diagram is irreducible if it is not reducible.
The length of a triangular diagram is the number of marked points on
the edge $XY$.
\section{Bijections}
We can now state the main theorem of this paper.
\begin{theorem}\label{mt} For all $n\ge 0$ there are inverse bijections
between words in $S$ of length $n$ and irreducible triangular diagrams of
length $n$.
\end{theorem}
Let $T(n)$ be the set of irreducible triangular diagrams of length
$n$ and $S^n$ the set of words in $S$ of length $n$. Then first we
construct a maps $T(n)\rightarrow S^n$ for $n\ge 0$. The construction
is essentially the same as the construction of the analogous map
in \cite{MR1680395}.

Choose a sequence of points in the edge $XY$, $(X_0,X_1,\ldots ,X_n)$, with
$X=X_0$ and $Y=X_n$ such that no point in the sequence is a marked
boundary point and such that for $1\le i\le n$ the interval $(X_{i-1},X_i)$
contains exactly one marked boundary point. Now for $0\le i\le n$ let
$\lambda_i$ be the weight of a minimal cut path with endpoints $A$ and
$X_i$. Then the claim is that $(\lambda_i-\lambda_{i-1})\in S$ for
$1\le i\le n$. This condition follows from the following:
\begin{prop}\label{td}
The irreducible triangular diagrams of length one
are precisely the seven triangular diagrams in
Figure \ref{bt}.
\end{prop}
An isoperimetric inequality for sectors is given in \cite{MR1800576}.
It would be interesting to know if this isoperimetric inequality
also holds for polyhedral metrics; and, if so, wether this implies
that the number of irreducible triangular diagrams of length one
is finite.

\newcommand{\TGC}{
\psset{linewidth=1pt,linestyle=dotted}
\pcline(-12,4)(12,4)
\pcline(12,4)(0,-8)\Aput{(1,0)}
\pcline(-12,4)(0,-8)\Bput{(0,2)}
\psset{linewidth=2pt,linestyle=solid}
\psline(-6,4)(-6,2)
\psline(-9,1)(-6,2)
\psline(-6,2)(0,0)
\psline(0,0)(6,2)
\psline(6,2)(9,1)
\psline(0,0)(0,-4)
\psline(0,-4)(3,-5)
\psline(0,-4)(-3,-5) }

\newcommand{\TGE}{
\psset{linewidth=1pt,linestyle=dotted}
\pcline(-12,4)(12,4)
\pcline(12,4)(0,-8)\Aput{(0,2)}
\pcline(-12,4)(0,-8)\Bput{(1,0)}
\psset{linewidth=2pt,linestyle=solid}
\psline(6,4)(6,2)
\psline(-9,1)(-6,2)
\psline(-6,2)(0,0)
\psline(0,0)(6,2)
\psline(6,2)(9,1)
\psline(0,0)(0,-4)
\psline(0,-4)(3,-5)
\psline(0,-4)(-3,-5) }

\newcommand{\TGG}{
\psset{linewidth=1pt,linestyle=dotted}
\pcline(-12,4)(12,4)
\pcline(12,4)(0,-8)\Aput{(0,0)}
\pcline(-12,4)(0,-8)\Bput{(0,1)}
\psset{linewidth=2pt,linestyle=solid}
\psline(0,4)(0,0)
\psline(0,0)(-6,-2) }

\newcommand{\TGD}{
\psset{linewidth=1pt,linestyle=dotted}
\pcline(-12,4)(12,4)
\pcline(12,4)(0,-8)\Aput{(0,1)}
\pcline(-12,4)(0,-8)\Bput{(0,1)}
\psset{linewidth=2pt,linestyle=solid}
\psline(0,4)(0,0)
\psline(0,0)(-6,-2)
\psline(0,0)(6,-2) }

\newcommand{\TGA}{
\psset{linewidth=1pt,linestyle=dotted}
\pcline(-12,4)(12,4)
\pcline(12,4)(0,-8)\Aput{(0,1)}
\pcline(-12,4)(0,-8)\Bput{(0,0)}
\psset{linewidth=2pt,linestyle=solid}
\psline(0,4)(0,0)
\psline(0,0)(6,-2) }

\newcommand{\TGF}{
\psset{linewidth=1pt,linestyle=dotted}
\pcline(-12,4)(12,4)
\pcline(12,4)(0,-8)\Aput{(0,1)}
\pcline(-12,4)(0,-8)\Bput{(1,0)}
\psset{linewidth=2pt,linestyle=solid}
\psline(-6,4)(-6,2)
\psline(0,0)(-6,2)
\psline(0,0)(0,-4)
\psline(3,-5)(0,-4)
\psline(-6,2)(-9,1)
\psline(0,-4)(-3,-5) }

\newcommand{\TGB}{
\psset{linewidth=1pt,linestyle=dotted}
\pcline(-12,4)(12,4)
\pcline(12,4)(0,-8)\Aput{(1,0)}
\pcline(-12,4)(0,-8)\Bput{(0,1)}
\psset{linewidth=2pt,linestyle=solid}
\psline(6,4)(6,2)
\psline(0,0)(6,2)
\psline(0,0)(0,-4)
\psline(-3,-5)(0,-4)
\psline(6,2)(9,1)
\psline(0,-4)(3,-5) }

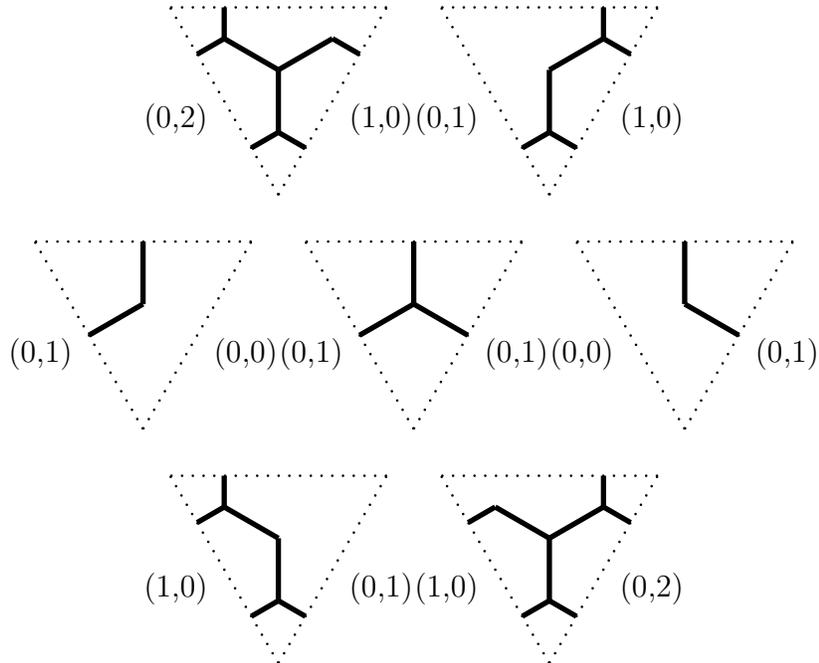
\begin{figure}
\begin{center}
\psset{xunit=0.1cm,yunit=0.1732cm}
\psset{xunit=1.2,yunit=1.2}
\begin{pspicture}(-45,-20)(45,20)
\rput(30,0){\TGA}
\rput(15,15){\TGB}
\rput(15,-15){\TGE}
\rput(-15,15){\TGC}
\rput(-15,-15){\TGF}
\rput(-30,0){\TGG}
\rput(0,0){\TGD}
\end{pspicture}
\psset{xunit=1cm,yunit=1cm}
\end{center}
\caption{Steps as triangles}\label{bt}
\end{figure}

\begin{proof} It is straightforward to check that each of the seven
diagrams in Figure \ref{bt} is an irreducible triangular diagram
of length one. It remains to show that these are all such diagrams.

First we observe that if we are given a triangular diagram with
one marked point on the top edge and such that the graph has
more than one component then the triangular diagram is reducible.

We assume that there is no region bounded by the graph.
This means we can now assume that we have a triangular diagram
with one marked point on the top edge and whose graph is a tree.
Then consider the graph obtained by removing all edges incident
to a marked point on the boundary. If this graph is not an 
interval then the original diagram can be reduced by pruning
the tree.

Note that if there are no marked points on the edge $AX$ then
the edge $AY$ has one marked point.
There is no diagram with one marked point on the top edge and no
other marked point. Also there cannot be more than one marked
point on the edge $AY$ otherwise the path $AXY$ would have weight
less than the weight of the path $AY$ and the cut path $AY$ is
required to be minimal. If there are two marked points on the
boundary then there is only one non-positive diagram.

Then to complete the proof it is sufficient to observe that
if you start with a diagram in Figure \ref{bt}, choose an
edge which does not meet the boundary, and add a new edge from
a point in this edge to one of the two sides of the triangle
then one of the following occurs:
\begin{enumerate}
\item The new diagram is already in Figure \ref{bt}.
\item In the new diagram either $AX$ or $AY$ is not a minimal
cut path.
\item The new triangular diagram is reducible.
\end{enumerate}

\end{proof}

Next we construct a map $S^n\rightarrow T(n)$.

Construct a planar graph by taking the subgraph of the square lattice 
on the vertices $(i,j)$ such that $0\le i,j$ and $i+j\le n$.
Then we label the edges of this graph by dominant weights.
Label the edge from $(i,j)$ to $(i,j+1)$ by $D_{i,j}$ and label
the edge from $(i,j)$ to $(i+1,j)$ by $H_{i,j}$. The labelling
is constructed by induction on $n-i-j$. If $i+j=n-1$ then the 
labels $D_{i,j}$ and $H_{i,j}$ are read off from the sequence
of elements of $S$.

If $D_{i+1,j-1}(k)\ge H_{i,j}(k)$ then put
\[ D_{i,j-1}(k)=D_{i+1,j-1}(k)- H_{i,j}(k)
\qquad
H_{i,j-1}=0 \]
and if $H_{i,j}(k)\ge D_{i+1,j-1}(k)$ then put
\[  D_{i,j-1}(k)0
\qquad
H_{i,j-1}=H_{i,j}(k)- D_{i+1,j-1}(k)
\]

\begin{figure}
\psset{xunit=0.10cm,yunit=0.1732cm}
\psset{xunit=0.75,yunit=0.75}
\begin{center}
\begin{pspicture}(-12,-20)(111,4)

\rput(0,0){
\psset{linewidth=1pt,linestyle=dotted}
\pcline(-12,4)(12,4)
\pcline(12,4)(0,-8)
\pcline(0,-8)(-12,4)\Aput{$D_A(i)$}

\pcline(12,4)(36,4)
\pcline(36,4)(24,-8)\Aput{$H_B(i)$}
\pcline(12,4)(24,-8)

\pcline(0,-8)(12,-20)\Bput{$D_B(i)-H_A(i)$}
\pcline(24,-8)(12,-20)\Aput{$0$}
}
\rput(0,0){$A$}
\rput(24,0){$B$}


\rput(75,0){
\psset{linewidth=1pt,linestyle=dotted}
\pcline(-12,4)(12,4)
\pcline(12,4)(0,-8)
\pcline(-12,4)(0,-8)\Bput{$D_A(i)$}

\pcline(12,4)(36,4)
\pcline(36,4)(24,-8)\Aput{$H_B(i)$}
\pcline(12,4)(24,-8)

\pcline(0,-8)(12,-20)\Bput{$0$}
\pcline(24,-8)(12,-20)\Aput{$H_A(i)-D_B(i)$}
\rput(0,0){$A$}
\rput(24,0){$B$}
}
\end{pspicture}
\end{center}
\psset{xunit=1cm,yunit=1cm}
\caption{Labelling the edges}\label{dw}
\end{figure}
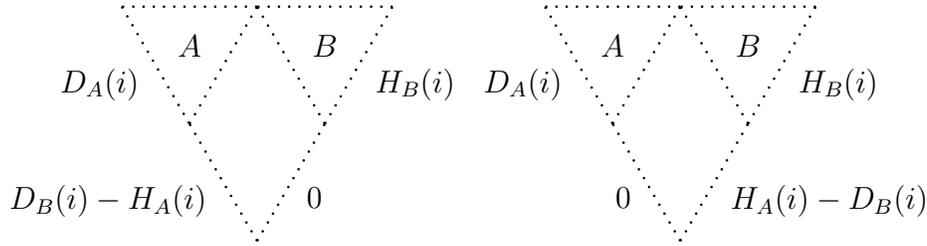

Now we come to draw the pictures. First we rotate so the
lines $i+j$ constant are horizontal. This means that the
labelling of the edges by dominant weights is given by the
two diagrams in Figure \ref{dw}.

Now we fill in the triangles and squares to make the diagram.
First we introduce a new type of edge. This will be drawn as a 
double edge. The Definition \ref{md} (of a diagram) is modified to allow
this new type of edge. The graph is still required to be trivalent but
we also allow a vertex to be incident to two of the original type of
edge and one of the new type. The Definition \ref{cp} (of a cut path) is
also modified to allow a cut path to intersect any edge transversally.
A cut path which crosses $a$ edges of the original type, contains
$b$ edges (of the the original type) and crosses $c$ edges of the
new type has weight $(a,b+c)$.

Then we modify the seven triangular diagrams in Figure \ref{bt}
and redraw them so that the edges are minimal cut paths which
do not contain any edge. This gives the seven triangular diagrams
in Figure \ref{nd}.
\newcommand{\TGU}{
\psset{linewidth=1pt,linestyle=dotted}
\pcline(-12,4)(12,4)
\pcline(12,4)(0,-8)\Aput{(1,0)}
\pcline(-12,4)(0,-8)\Bput{(0,1)}
\psset{linewidth=2pt,linestyle=solid}
\psline(0,4)(0,0)
\psline(0,0)(-6,-2)
\psline[doubleline=true,linewidth=1pt](0,0)(6,-2) }

\newcommand{\TGV}{
\psset{linewidth=1pt,linestyle=dotted}
\pcline(-12,4)(12,4)
\pcline(12,4)(0,-8)\Aput{(0,1)}
\pcline(-12,4)(0,-8)\Bput{(1,0)}
\psset{linewidth=2pt,linestyle=solid}
\psline(0,4)(0,0)
\psline[doubleline=true,linewidth=1pt](0,0)(-6,-2)
\psline(0,0)(6,-2) }

\newcommand{\TGX}{
\psset{linewidth=1pt,linestyle=dotted}
\pcline(-12,4)(12,4)
\pcline(12,4)(0,-8)\Aput{(1,0)}
\pcline(-12,4)(0,-8)\Bput{(0,2)}
\psset{linewidth=2pt,linestyle=solid}
\psline(-6,4)(-6,2)
\psline(0,0)(-6,2)
\psline(0,0)(0,-4)
\psline(-3,-5)(0,-4)
\psline(-6,2)(-9,1)
\psline[doubleline=true,linewidth=1pt](0,-4)(3,-5) }

\newcommand{\TGY}{
\psset{linewidth=1pt,linestyle=dotted}
\pcline(-12,4)(12,4)
\pcline(12,4)(0,-8)\Aput{(0,2)}
\pcline(-12,4)(0,-8)\Bput{(1,0)}
\psset{linewidth=2pt,linestyle=solid}
\psline(6,4)(6,2)
\psline(0,0)(6,2)
\psline(0,0)(0,-4)
\psline(3,-5)(0,-4)
\psline(6,2)(9,1)
\psline[doubleline=true,linewidth=1pt](0,-4)(-3,-5) }

\begin{figure}
\begin{center}
\psset{xunit=0.1cm,yunit=0.1732cm}
\psset{xunit=1.2,yunit=1.2}
\begin{pspicture}(-45,-20)(45,20)
\rput(30,0){\TGA}
\rput(15,-15){\TGY}
\rput(15,15){\TGU}
\rput(-15,-15){\TGV}
\rput(-15,15){\TGX}
\rput(-30,0){\TGG}
\rput(0,0){\TGD}
\end{pspicture}
\psset{xunit=1cm,yunit=1cm}
\end{center}
\caption{Steps as triangles}\label{nd}
\end{figure}
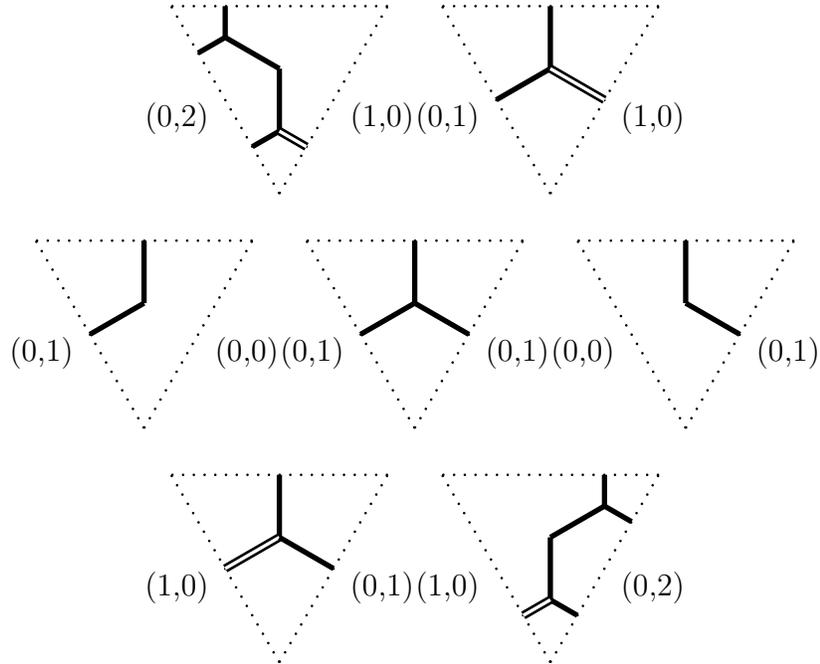

Then given a word in $S$ we first draw the grid as above.
Then we fill in each triangle using the previous diagram.
Next we fill in the diamonds. Note that the weights $H_A$
and $D_B$ are both elements of the set 
\[ \{(0,0),(0,1),(0,2),(1,0)\} \]
This gives sixteen different diamonds. Furthermore note
that for each of these sixteen diamonds the other two
weights are also elements of this set. This implies that
for any word every diamond will have one of these sixteen
labellings. Therefore to complete the diagram it is sufficient
to know how to fill in a diamond with each of these labellings.

There are four symmetric diamonds which are also the diamonds
where the weights on the top two edges are equal and the weights
on the bottom two edges are zero. The case in which all four weights
are zero is an empty diamond. The other three diamonds are

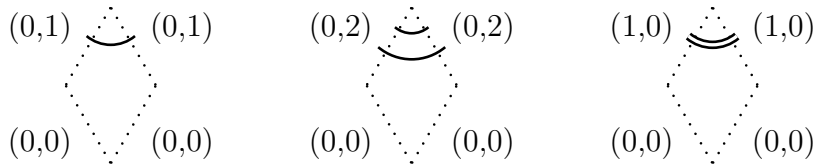
\begin{center}
\psset{xunit=0.1cm,yunit=0.1732cm}
\begin{pspicture}(-6,-8)(86,4)
\rput(0,0){
\psset{linewidth=1pt,linestyle=dotted}
\pcline(0,4)(6,-2)\Aput{(0,1)}
\pcline(6,-2)(0,-8)\Aput{(0,0)}
\pcline(0,-8)(-6,-2)\Aput{(0,0)}
\pcline(-6,-2)(0,4)\Aput{(0,1)}
\psset{linestyle=solid}
\psarc(0,4){0.5}{-130}{-50}
}
\rput(40,0){
\psset{linewidth=1pt,linestyle=dotted}
\pcline(0,4)(6,-2)\Aput{(0,2)}
\pcline(6,-2)(0,-8)\Aput{(0,0)}
\pcline(0,-8)(-6,-2)\Aput{(0,0)}
\pcline(-6,-2)(0,4)\Aput{(0,2)}
\psset{linestyle=solid}
\psarc(0,4){0.3464}{-130}{-50}
\psarc(0,4){0.6928}{-130}{-50}
}
\rput(80,0){
\psset{linewidth=1pt,linestyle=dotted}
\pcline(0,4)(6,-2)\Aput{(1,0)}
\pcline(6,-2)(0,-8)\Aput{(0,0)}
\pcline(0,-8)(-6,-2)\Aput{(0,0)}
\pcline(-6,-2)(0,4)\Aput{(1,0)}
\psset{linestyle=solid}
\psarc[doubleline=true,linewidth=1pt](0,4){0.5}{-130}{-50}
}
\end{pspicture}
\psset{xunit=1cm,yunit=1cm}
\end{center}
The other twelve diamonds come in pairs and we only give one
member of each pair. The other member is obtained by reflection
in a vertical line. There are three diagrams in which an oposite
pair of edges have weight zero. These give the three diamonds

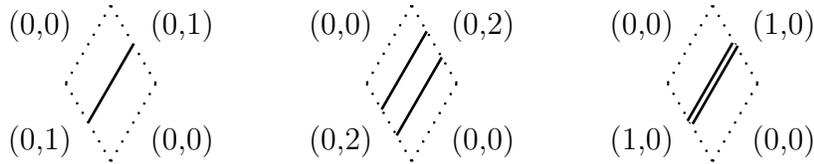
\begin{center}
\psset{xunit=0.1cm,yunit=0.1732cm}
\begin{pspicture}(-6,-8)(86,4)
\rput(0,0){
\psset{linewidth=1pt,linestyle=dotted}
\pcline(0,4)(6,-2)\Aput{(0,1)}
\pcline(6,-2)(0,-8)\Aput{(0,0)}
\pcline(0,-8)(-6,-2)\Aput{(0,1)}
\pcline(-6,-2)(0,4)\Aput{(0,0)}
\psset{linestyle=solid}
\psline(3,1)(-3,-5)
}
\rput(40,0){
\psset{linewidth=1pt,linestyle=dotted}
\pcline(0,4)(6,-2)\Aput{(0,2)}
\pcline(6,-2)(0,-8)\Aput{(0,0)}
\pcline(0,-8)(-6,-2)\Aput{(0,2)}
\pcline(-6,-2)(0,4)\Aput{(0,0)}
\psset{linestyle=solid}
\psline(2,2)(-4,-4)
\psline(4,0)(-2,-6)
}
\rput(80,0){
\psset{linewidth=1pt,linestyle=dotted}
\pcline(0,4)(6,-2)\Aput{(1,0)}
\pcline(6,-2)(0,-8)\Aput{(0,0)}
\pcline(0,-8)(-6,-2)\Aput{(1,0)}
\pcline(-6,-2)(0,4)\Aput{(0,0)}
\psset{linestyle=solid}
\psline[doubleline=true,linewidth=1pt](3,1)(-3,-5)
}
\end{pspicture}
\psset{xunit=1cm,yunit=1cm}
\end{center}

The remaining three diamonds are
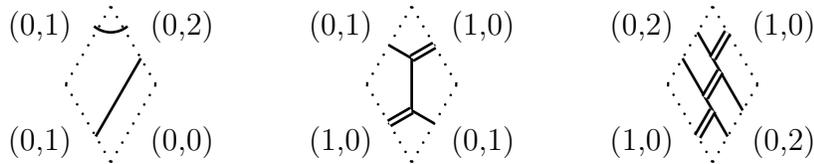
\begin{center}
\psset{xunit=0.1cm,yunit=0.1732cm}
\begin{pspicture}(-6,-8)(86,4)

\rput(0,0){
\psset{linewidth=1pt,linestyle=dotted}
\pcline(0,4)(6,-2)\Aput{(0,2)}
\pcline(6,-2)(0,-8)\Aput{(0,0)}
\pcline(0,-8)(-6,-2)\Aput{(0,1)}
\pcline(-6,-2)(0,4)\Aput{(0,1)}
\psset{linestyle=solid}
\psarc(0,4){0.3464}{-130}{-50}
\psline(4,0)(-2,-6)
}
\rput(40,0){
\psset{linewidth=1pt,linestyle=dotted}
\pcline(0,4)(6,-2)\Aput{(1,0)}
\pcline(6,-2)(0,-8)\Aput{(0,1)}
\pcline(0,-8)(-6,-2)\Aput{(1,0)}
\pcline(-6,-2)(0,4)\Aput{(0,1)}
\psset{linestyle=solid}
\psline(-3,1)(0,0)
\psline[doubleline=true](3,1)(0,0)
\psline(0,0)(0,-4)
\psline[doubleline=true](0,-4)(-3,-5)
\psline(0,-4)(3,-5)
}
\rput(80,0){
\psset{linewidth=1pt,linestyle=dotted}
\pcline(0,4)(6,-2)\Aput{(1,0)}
\pcline(6,-2)(0,-8)\Aput{(0,2)}
\pcline(0,-8)(-6,-2)\Aput{(1,0)}
\pcline(-6,-2)(0,4)\Aput{(0,2)}
\psset{linestyle=solid}
\psline(-2,2)(4,-4)
\psline(-4,0)(2,-6)
\psline[doubleline=true](2,2)(0,0)
\psline[doubleline=true](1,-1)(-1,-3)
\psline[doubleline=true](0,-4)(-2,-6)
}
\end{pspicture}
\psset{xunit=1cm,yunit=1cm}
\end{center}

In the resulting diagram each interior vertex is still trivalent
and each interior trivalent vertex is either incident to three
of the original type of edge or two of the original type of edge
and one edge of the new type. Hence each ocurrence of the new
type of edge can be removed by the following replacement.

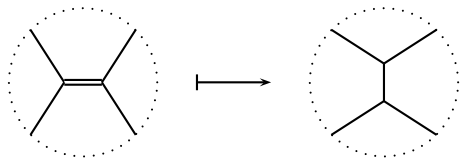
\begin{center}
\begin{pspicture}(4,-1)(11,2)
\psset{linestyle=dotted}
\pscircle(6,0){1}
\pscircle(10,0){1}
\psset{linestyle=solid}
\psline(5.3,0.7)(5.75,0)
\psline[doubleline=true](5.75,0)(6.25,0)
\psline(6.25,0)(6.7,-0.7)
\psline(5.75,0)(5.3,-0.7)
\psline(6.7,0.7)(6.25,0)

\psline(9.3,0.7)(10,0.25)
\psline(10,0.25)(10,-0.25)
\psline(10,-0.25)(9.3,-0.7)
\psline(10.7,0.7)(10,0.25)
\psline(10,-0.25)(10.7,-0.7)
\psline[arrows=|->](7.5,0)(8.5,0)
\end{pspicture}
\end{center}

The reason we have to introduce the new type of edge and the
modified pictures for the elements of $S$ is that if we used the
original pictures then we could get a diagram with squares.
When working by hand it is straightforward to remove these 
superfluous squares. Informally, the rule is that whenever we see a ladder
with more than one rung, to remove all but one rung.

An important feature of this construction is the following observation:
\begin{lemma}\label{cl} Any increasing path which
starts and the bottom vertex and follows the construction lines to the top
edge of the diagram is a minimal cut path.
\end{lemma}

This observation shows that the diagram is non-positive since it shows
that we have a minimal cut path from the bottom vertex to the top edge
which passes through any internal region. Hence we have constructed
a triangular diagram of length $n$. It remains to check that this is
irreducible. This is the statement that the edges of the triangle
are the unique minimal cut paths between their endpoints.

Then Lemma \ref{cl} also shows that if we start with a word
construct the diagram and then derive a word that we recover the 
original word. In particular this shows that  for $n\ge 0$, the
map $S^n\rightarrow T(n)$ is surjective.

Next we show that, for $n\ge 0$, the map $S^n\rightarrow T(n)$ is injective.
The proof is by induction on the length of the word. The basis
of the induction is the case of length one which is proved in
Proposition \ref{td}. Assume the result for words of length $n$.
Let $w$ be a word of length $n+1$. Let $w_i$ be obtained by dropping
the final step and $w_f$ be obtained by dropping the first step.
Then these are both words oflength $n$ and so by the inductive
hypothesis have unique diagrams.
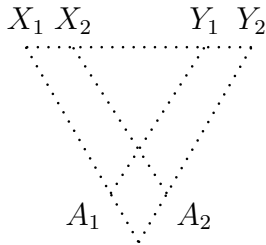
\begin{figure}
\psset{xunit=0.10cm,yunit=0.1732cm}
\psset{xunit=1.25,yunit=1.25}
\begin{center}
\begin{pspicture}(-12,-8)(12,6)

\psset{linewidth=1pt,linestyle=dotted}
\pcline(-12,4)(12,4)
\pcline(12,4)(0,-8)
\pcline(0,-8)(-12,4)

\pcline(-7,4)(3,-5)
\pcline(7,4)(-3,-5)
\uput[90](-12,4){$X_1$}
\uput[90](12,4){$Y_2$}
\uput[90](-7,4){$X_2$}
\uput[90](7,4){$Y_1$}
\uput[225](-3,-5){$A_1$}
\uput[-45](3,-5){$A_2$}
\end{pspicture}
\end{center}
\psset{xunit=1cm,yunit=1cm}
\caption{Inductive step}\label{iw}
\end{figure}
This means that in Figure \ref{iw} the triangle $A_1X_1Y_1$
can be filled in uniquely using the word $w_i$ and the triangle
$A_2X_2Y_2$ can be filled in uniquely using the word $w_f$.
Thus the only part of the diagram that has not been filled in
is the lowest diamond between the points $A_1$ and $A_2$.
Thus the claim is that for each of the sixteen possible diamonds
the rule we have given for filling it in is the unique rule
that gives a non-positive irreducible triangular diagram.
For each of these sixteen diamonds, the claim can be checked
using \cite[Lemma 6.5]{MR1403861}.

This proves Theorem \ref{mt}. Then for $n\ge 0$, these bijections
restrict to bijections between the set of non-positive diagrams
with $n$ boundary points and the set of lattice walks of length $n$
which start and end at the origin.

\bibliographystyle{halpha}
\bibliography{g2bij.bib}

\end{document}